\documentclass[a4paper]{amsart}
\usepackage{amssymb}
\usepackage{amsmath}
\usepackage{amsfonts}
\usepackage{graphicx}
\usepackage{pstricks}
\theoremstyle{plain}
\newtheorem{theo}{Theorem}[section]
\newtheorem{conj}[theo]{Conjecture}

\def\diam{\mathop{\rm diam}\nolimits}
\begin{document}
                  
\bigskip

\title[On near-embeddings]{Constructing near-embeddings of codimension one manifolds with countable dense singular sets}

\author{D. Repov\v s}
\address{Institute of Mathematics, Physics and Mechanics,
and Faculty of Education, University of Ljubljana, P.O.Box 2964,
Ljubljana 1001, Slovenia}
\email{dusan.repovs@guest.arnes.si}

\author{W. Rosicki}
\address{Institute of Mathematics, Gdansk University, ul. Wita Stwosza 57, 80-952 Gda\'nsk, Poland}
\email{wrosicki@math.univ.gda.pl}

\author{A. Zastrow}
\address{Institute of Mathematics, Gdansk University, ul. Wita Stwosza 57, 80-952 Gda\'nsk, Poland}
\email{zastrow@math.univ.gda.pl}

\author{M. \v Zeljko}
\address{Institute of Mathematics, Physics and Mechanics,
University of Ljubljana, P.O.Box 2964,
Ljubljana 1001, Slovenia}
\email{matjaz.zeljko@fmf.uni-lj.si}

\date{March 29, 2008}



\begin{abstract}
The purpose of this paper is to present, for all $n\ge 3$,  very simple examples of continuous maps $f:M^{n-1} \to M^{n}$ from  closed $(n-1)$-manifolds $M^{n-1}$
into
 closed
$n$-manifold $M^n$
such that even though
the singular set $S(f)$ of $f$
is countable and dense,  the map $f$ can nevertheless be approximated by an embedding, i.e. $f$ is a {\sl near-embedding}.  Notice that in such case in dimension 3 one can get even a piecewise-linear approximation by an embedding \cite{Bing}.
\end{abstract}

\maketitle

\section{Introduction}
Denote the {\em singular} set of an arbitrary continuous mapping 
$f\colon X\to Y$ between topological spaces
by 
$S(f)=\{x\in X \mid f^{-1}(f(x))\ne x\}$. 
A manifold which is connected, compact and has no boundary is said to be {\sl closed}.
The following conjecture was proposed by the first author in the mid 1980's. 
For $n=3$ it has since been shown to be equivalent to the Bing Conjecture from the 1950's  (cf.  \cite{Brahm1}) and it
is closely related to the 3-dimensional Recognition Problem, one of the central problems of geometric topology (cf. \cite{Repovs}). In particular, it is closely related to the general position of 3-manifolds, called the {\sl  light map separation property LMSP$^*$}, introduced by Daverman and Repov\v{s} 
(cf. Conjecture 5.4 in \cite{DaRe}):

\begin{conj}
Let $f\colon M^{n-1} \to M^n$ be any continuous (posibly surjective) map from
a closed $(n-1)$-manifold $M^{n-1}$
into a closed n-manifold $M^n$, $n\ge 3$,
such that $\dim S(f)=0$. 
Then for every $\varepsilon >0$ there exists an embedding $g\colon M^{n-1}\to M^n$ such that for every $x\in M^{n-1}$, 
$d(f(x),g(x))<\varepsilon$, i.e. $f$ is a near-embedding.
\end{conj}
 
In the case when $n=3$, a very special case of Conjecture 1.1 
was verified by Anderson \cite{Anderson} in 1965. Then in 1992
Brahm \cite{Brahm2} proved  Conjecture 1.1 
for the case when $n=3$ and the closure of the singular set is 0-dimensional, dim(Cl$S(f))=0$.   
However, in general the second  property needs not be satisfied -- 
it was shown in  \cite{Brahm3} that for $n=3$ it can happen  that dim(Cl$S(f))=n-1$.

Since the construction in \cite{Brahm3} is very technical, there has been for a long time an open question if there is an {\sl elementary} example of a continuous map $f:M^2\to M^3$ such that $\dim S(f)=0$ whereas $0<\dim$Cl$(S(f))\le 3$. 
The purpose of this note is to present such an example -- it has a very simple construction and the verification of all  asserted properties is straightforward. 
Moreover, unlike \cite{Brahm3}, our methods evidently  generalize in a direct manner
to yield continuous maps $f\colon M^{n-1} \to M^{n}$ of 
closed codimension one manifolds into  closed $n$-manifolds, with properties analogous to (i) and (ii) below, for every $n\ge 3$. 

\begin{theo}
For every $n\ge 3$, there exists a continuous map $f\colon S^{n-1}\to S^n$ such that:

{\rm{(i)}} the singular set
$S(f)$ of $f$ is countable and dense  (hence 0-dimensional and nonclosed); and

 {\rm{(ii)}}$f$ can be approximated by an embedding. 
\end{theo} 

{\sl Remark.} A question when a light map is a near-embedding is also interesting in view of the classical {\sl Monotone-Light Factorization Theorem} (cf. e.g. \cite{Whyburn}) which asserts that
that every continuous mapping 
$f:X \to Y$ from any compact space $X$ to any space $Y$ can be factorized as a
product $f=l \circ m$ 
of a monotone map $m:X\to Z$ (i.e. each point  inverse $m^{-1}(z)$ is connected)
and 
a  light map $l:Z \to Y$ (i.e. each point inverse $l^{-1}(y)$   is totally disconnected). 


\section{Proof of Theorem 1.2}

Let $n\ge 3$ and choose   a
countable basis $\{ U_i\}_{i\in N}$
of open sets for
$S^{n-1}$ 
(which is considered as the standardly embedded $(n-1)$-sphere in $S^{n}$). 
We shall inductively construct a sequence of pairwise disjoint 
{\sl tame} PL arcs $\alpha_i$ in $S^{n}$ 
(i.e.  for every $i\ge 1$ there is a homeomorphism 
$h_i:S^{n} \to S^{n}$ such that $h_i(\alpha_i) \subset S^{1} \subset S^{n}$) 
with the property 
that:
\begin{enumerate}
\item
for every $i$, $\alpha_i\cap S^{{n-1}}=\partial \alpha_i\subset U_i$; 
\item
for every $i$, $\diam(\alpha_i)<\frac{1}{2^i}$.
\end{enumerate}

Begin with a tame PL arc  $\alpha_1 \subset S^{n}$ such that
$\partial \alpha_1 \subset U_1$ and $\diam (\alpha_1)<1/2$. Assume inductively,  that we have 
already constructed pairwise disjoint tame PL arcs
$\alpha_1$, \ldots, $\alpha_{n-1}\subset S^{n}$ with all required properties. 
We can then clearly find a tame PL arc $\alpha_n\subset S^{n}$ such that 
$\partial \alpha_n\subset U_n$ and $\alpha_n$ is disjoint with $\alpha_1\cup\ldots\cup \alpha_{n-1}$.

\begin{figure}[ht]
\center
\includegraphics{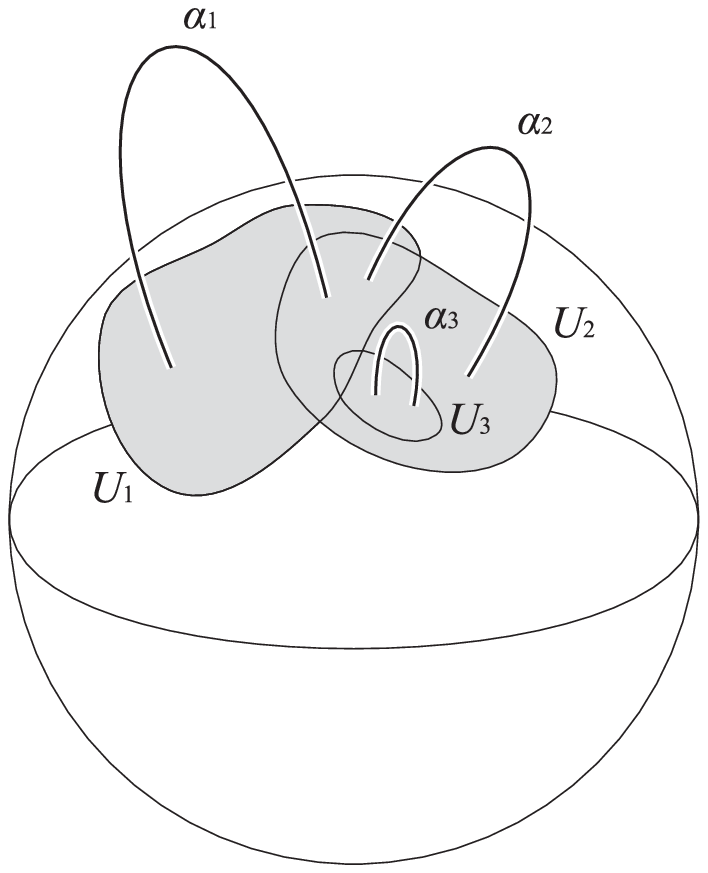}
\end{figure}

By our construction, $\{\alpha_n\}_{n\in N}$ is a {\sl null-sequence}, i.e. lim$_{i \to \infty}$$\diam{\alpha_{i}}=0$. 
The decomposition ${{G}}=\{\alpha_i\}_{n\in N}$ of $S^{n}$ into points and arcs is 
clearly {\sl cellular} (i.e. each element of the decomposition ${G}$
is the intersection of a {\sl nested} sequence of closed $n$-cells $\{B^{n}_{k}\}_{k\in N}$
in $S^{n}$) (i.e. for ever $k, B^{n}_{k+1}\subset$ Int$B^{n}_{k}$)
and {\sl upper semicontinuous} (i.e. the quotient map 
$\pi:S^{n} \to S^{n}/{{G}}$ is closed). 

Therefore it follows by \cite[Theorem 7, page 56]{Daverman} that the decomposition $G$ is {\sl shrinkable} (i.e. the map $\pi$ is approximable by homeomorhisms). In particular, the quotient space $S^{n}/G$ is homeomorphic to $S^{n}$.
The desired
mapping 
$f\colon S^{n-1}\to S^{n}$
is now defined as 
the compositum 
$f= \pi \circ i$
of the inclusion $i\colon S^{n-1}\to S^{n}$ and  
the decomposition
quotient mapping $\pi\colon S^{n}\to S^{n}/G$. 

It follows by 
construction that 
the singular set
$S(f)=\bigcup_i \partial \alpha_i$
is 
countable and dense. It is also clear that this map can be approximated arbitrarily closely by embeddings of $S^{n-1}$ into $S^{n}$ (by not shrinking the arcs all the way but only as much as it is necessary to make them sufficiently small).

\section{Acknowledgements}
This research was supported by the Polish--Slovenian  grant BI-PL 10/2008-2009. The first and the fourth author were supported by the ARRS program P1-0292-0101-04 and project J1-9643-0101. The second and the third
authors were partially supported by KBN grant N200100831/0524.

\end{document}